# A new zero-order 1-D optimization algorithm: trichotomy method

Alena Antonova
Lyceum 77
Chelyabinsk, Russia
a_antonova02@mail.ru

Olga Ibryaeva
School of Electrical Engineering and Computer Science
South Ural State University
Chelyabinsk, Russia
ibriaevaol@susu.ru

**Abstract:** A new 1D search method is proposed for minimizing an arbitrary real valued function. The algorithm is a modification of the interval halving method which is based on dividing the interval of uncertainty by three points into four equal parts. The trichotomy method is based on dividing the interval by five points into six equal parts and provides the interval reducing exactly three times in every iteration. New algorithm's performance has been extensively tested and compared to well-known 1D search algorithms: the interval halving method, golden section method, Fibonacci search. The results show that the trichotomy method usually require less calculations of the values of the minimized function to determine the minimum point with the given accuracy. Also it has better accuracy when finding the minimum point after using $N$ calculations of the function values.

## I. Introduction

We are given a function $f(x): R \to R$ and an interval $[a, b]$ in which this function is known to have one minimum. We are to locate the minimum with some specified accuracy.
Consider the following optimization problem:
$$\min_{x \in R} f(x).$$
Let $x^* \in [a, b]$ be the minimum point of $f(x)$.
The function $f(x)$ is said to be *unimodal* on $[a, b]$, i.e. for $a \le x_1 < x_2 \le b$:
$$x_2 < x^* \quad \Rightarrow \quad f(x_1) > f(x_2),$$
$$x_1 > x^* \quad \Rightarrow \quad f(x_2) > f(x_1).$$
Some examples of unimodal functions are shown in Figs. 1-2. Thus a unimodal function can be a nondifferentiable or even a discontinuous function.
The dichotomous search method, as well as the Fibonacci and the golden section methods, are sequential search methods in which the result of any experiment influences the location of the subsequent experiment. In the dichotomous search, two experiments are placed as close as possible at the center of the interval of uncertainty. Based on the relative values of the objective function at the two points, almost half of the interval of uncertainty is eliminated.

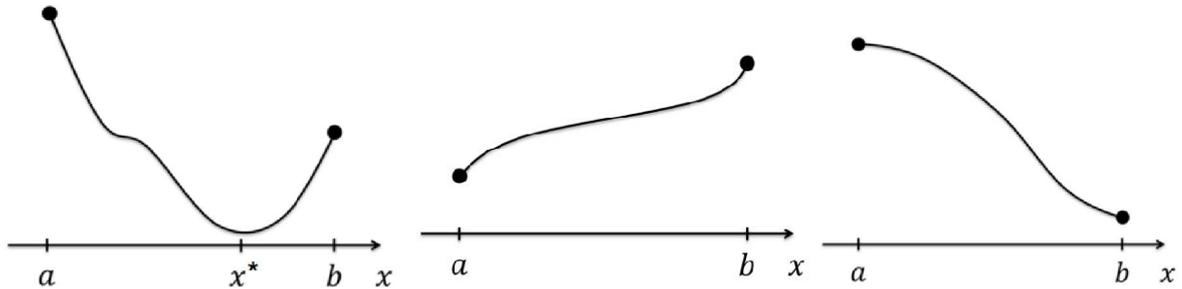
Figure 1. Unimodal functions

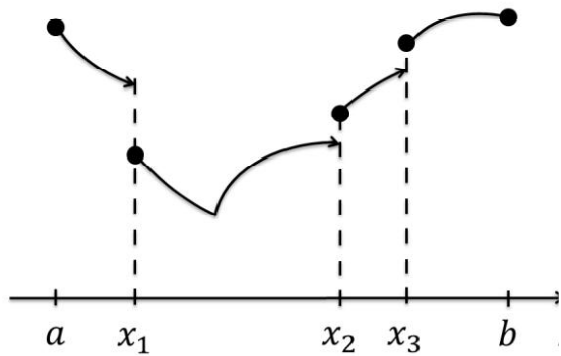
Figure 2. Discontinuous unimodal function

We will not discuss dichotomous search in detail because our new method is actually a generalization of a similar method – Interval halving method described in next section. A complete description of dichotomous search can be found, for example, in [1].

The rest of this paper is organized as follows. The interval halving method is described in Section 2. Our new algorithm is introduced in Section 3. Then, in Section 4, we obtain some theoretical estimates of the effectiveness of the proposed method in comparison with the interval halving method. In Section 5, an experimental study of the algorithm's performance is presented. The conclusions are presented in Section 6.

## II. Interval halving method

In the interval halving method, exactly one-half of the current interval of uncertainty is deleted in every stage. It requires at least three experiments in the first stage and at least two experiments in each subsequent stage. The procedure can be described by the following steps:

**1.** Find $x_2 = \frac{a+b}{2}$, $f(x_2)$, and go to step 2.
**2.** Find $x_1 = \frac{a+x_2}{2}$, $f(x_1)$, and go to step 3.
**3.** Compare $f(x_1)$ with $f(x_2)$.

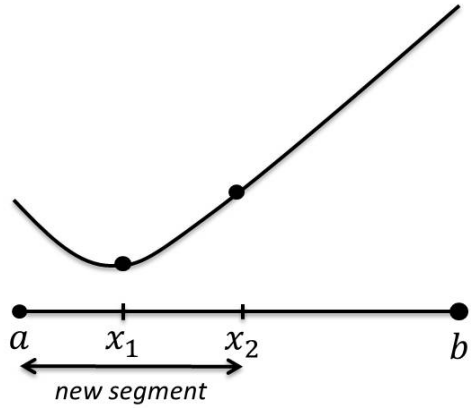

If $f(x_1) \le f(x_2)$, as shown in Fig. 3, then the new interval of uncertainty will be $[a; x_2]$. Set $b = x_2$, $x_2 = x_1$, $f(x_2) = f(x_1)$, and go to step 5.

Else – find $x_3 = \frac{x_2+b}{2}$, $f(x_3)$, and go to step 4.

**Figure 3.** $f(x_1) \le f(x_2)$

**4.** Compare $f(x_2)$ with $f(x_3)$.
If $f(x_2) \le f(x_3)$, as shown in Fig. 4, then the new interval of uncertainty will be $[x_1; x_3]$. Set $a = x_1$, $b = x_3$.

Else – the new interval of uncertainty will be $[x_2; b]$ (Fig. 5). Set $a = x_2$, $x_2 = x_3$, $f(x_2) = f(x_3)$.
Go to step 5.

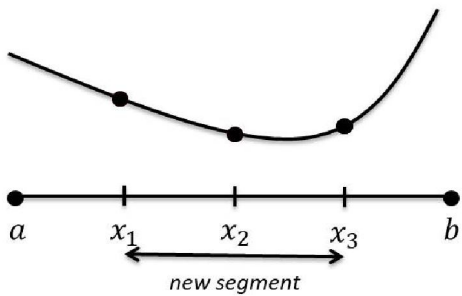
**Figure 4.** $f(x_2) \le f(x_3)$

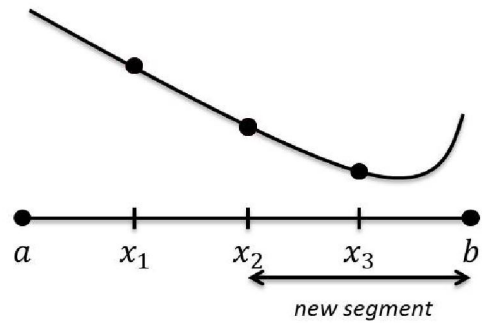
**Figure 5.** $f(x_2) > f(x_3)$

**5.** Test whether the convergence criterion is satisfied.
Find $\varepsilon_n = \frac{b-a}{2}$ and compare with $\varepsilon$. If $\varepsilon_n > \varepsilon$, then go to step 2. Else – stop the procedure, set $x^* \approx x_2$, $f^* \approx f(x_2)$.

The above method is based on dividing an interval by three points into four equal parts. In the next section, we will propose its modification in the case of dividing an interval by five points into six equal parts.

## III. Trichotomy method

Start, as in the interval halving method, from the middle of the interval.

**1.** Find $x_3 = \frac{a+b}{2}$, $f(x_3)$ and go to step 2.
**2.** Find $x_2 = \frac{a+2x_3}{3}$, $f(x_2)$ and go to step 3.

**3.** Compare $f(x_2)$ with $f(x_3)$.

If $f(x_2) \leq f(x_3)$, then set $x_1 = \frac{a+x_2}{2}$, and find $f(x_1)$.

    If $f(x_1) \leq f(x_2)$, as shown in Fig. 6, then the new interval of uncertainty will be $[a; x_2]$. Set $b = x_2, x_3 = x_1, f(x_3) = f(x_1)$.

    Else – the new interval of uncertainty will be $[x_1; x_3]$. Set $a = x_1, b = x_3, x_3 = x_2, f(x_3) = f(x_2)$ (Fig. 7).

Go to step 5.

Else (if $f(x_2) > f(x_3)$) – find $x_4 = \frac{b+2x_3}{3}, f(x_4)$, and go to step 4.

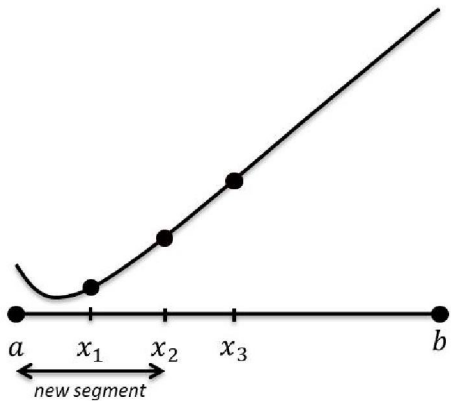 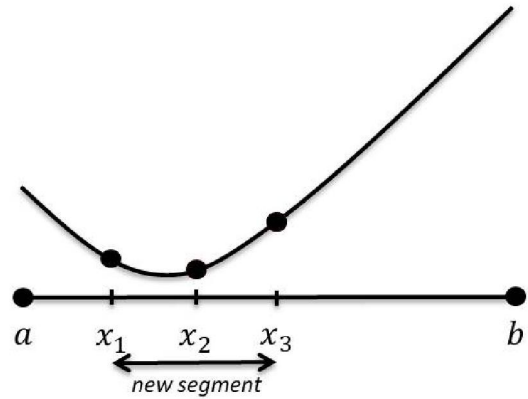

**Figure 6.** $f(x_1) \leq f(x_2)$      **Figure 7.** $f(x_1) > f(x_2)$

**4.** Compare $f(x_3)$ with $f(x_4)$.

If $f(x_4) \leq f(x_3)$, find $x_5 = \frac{2b+x_3}{3}, f(x_5)$.

    If $f(x_5) \leq f(x_4)$, as shown in Fig. 8, then the new interval of uncertainty will be $[x_4; b]$. Set $a = x_4, x_3 = x_5, f(x_3) = f(x_5)$.

    Else – the new interval of uncertainty will be $[x_3; x_5]$. Set $a = x_3, b = x_5, x_3 = x_4, f(x_3) = f(x_4)$ (Fig. 9).

Go to step 5.

Else (if $f(x_4) > f(x_3)$), as shown in Fig. 10, then the new interval of uncertainty will be $[x_2; x_4]$. Set $a = x_2, b = x_4$. Go to step 5.

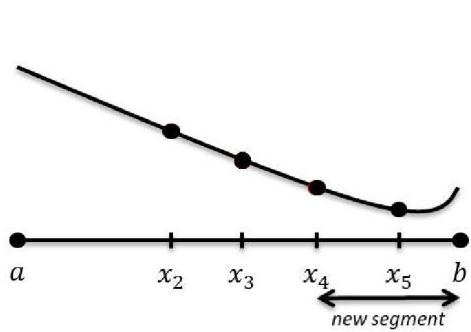 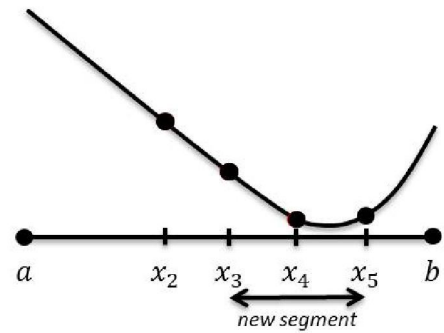

**Figure 8.** $f(x_5) \leq f(x_4)$      **Figure 9.** $f(x_5) > f(x_4)$

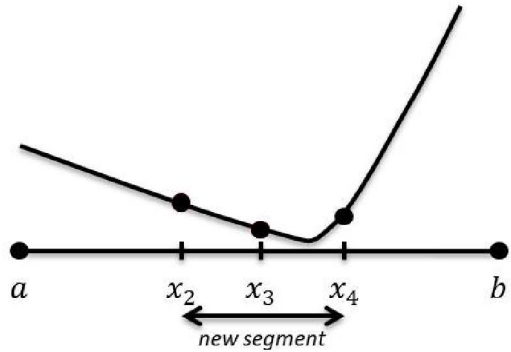

**Figure 10.** $f(x_4) > f(x_3)$

5. Test whether the convergence criterion is satisfied.

Find $\varepsilon_n = \frac{b-a}{2}$ and compare with $\varepsilon$. If $\varepsilon_n > \varepsilon$, then go to step 2. Else – stop the procedure, set $x^* \approx x_3$, $f^* \approx f(x_3)$.

In the method above, the interval of uncertainty is reduced exactly three times in every stage. It requires no more than five calculations of the function values in the first stage and no more than four calculations in each subsequent stage.

Initially, taking into account that function evaluation is usually an expensive operation, it seems that the proposed method will yield to interval halving method, since will require the calculation of more function values. However, in reality, it is rarely required to calculate exactly four function values at each iteration and, as we will see further from the results of numerical experiments, the trichotomy method provides greater accuracy in finding the minimum for a fixed number $N$ of calculations of the function values, rather than the interval halving method. This is also due to the faster reduction of the interval of uncertainty (in three, not two times) by the proposed method.

## IV. Theoretical estimation of the effectiveness of the interval halving and trichotomy methods

In this section, we obtain some theoretical estimates of the effectiveness of the proposed method in comparison with the interval halving method. The following paragraph is devoted to the results of numerical experiments.

The following criteria were chosen as criteria for evaluating the effectiveness of the algorithms:

1. The number of iterations of a method required to determine the minimum point with a given accuracy. A more effective method is one that requires fewer steps.

2. The accuracy with which we find the minimum point after using N calculations of the values of the minimized function. It should be noted that often in applications, the calculation of a function value is associated with expensive physical experiments, or the function itself has a complex form, which makes its calculation time-consuming. From this point of view, the best algorithm is the one that gives the same accuracy with fewer requests to the function calculation.

Propositions 1 and 2 compare two methods based on these criteria.

**Proposition 1.** *The number of iterations of the interval halving method $k_D$ and the trichotomy method $k_T$, which is necessary for determining the point $x^*$ with the accuracy $\varepsilon$, satisfies the inequalities:*

$$k_D \geq \left[\log_2 \frac{b-a}{2\varepsilon}\right] + 1, \qquad k_T \geq \left[\log_3 \frac{b-a}{2\varepsilon}\right] + 1,$$

*respectively. Here brackets [...] mean the integer part of a number.*

**Proof.** At each iteration of the interval halving method, the length of the interval of uncertainty is halved, that is, after $k$ iterations we have an interval of length $\frac{b-a}{2^k}$. Then $\frac{b-a}{2^k} \leq \varepsilon$, whence $2^k \geq \frac{b-a}{2\varepsilon}$ and $k \geq \log_2 \frac{b-a}{2\varepsilon}$. Since the number of iterations must be integer, we come to the formula $k_D \geq \left[\log_2 \frac{b-a}{2\varepsilon}\right] + 1$.

A similar formula for the method of trichotomy follows from the reduction of the interval of uncertainty three times at each step of this method.

**Note 1:** As one can see, the trichotomy method requires less iteration than the interval halving method. This result is quite expected, since at each step of the trichotomy method, the current interval of uncertainty is reduced three times, and in the interval halving method, only twice.

**Proposition 2.** *The accuracy $\varepsilon$ with which we find the minimum point $x^*$ after using $N$ calculations of the values of the minimized function $f(x)$ can be estimated as follows:*

$$\varepsilon_D \leq \frac{b-a}{2 \cdot 2^{\frac{N-1}{2}}}, \qquad \varepsilon_T \leq \frac{b-a}{2 \cdot 3^{\frac{N-1}{4}}}.$$

*for the interval halving method and the trichotomy method, respectively.*

**Proof.** For the interval halving method, no more than three values of the function $f(x)$ are calculated in the first step and no more than two values – in the following steps.

Thus, $k$ iterations require
$$N \leq 3 + 2(k-1) = 2k + 1$$
calculations of the $f(x)$.

Since the minimum number of calculations of the function values is two in the first step and one in $k-1$ subsequent steps, then $N \geq 2 + k - 1 = k + 1$. So, using $N$ calculations of the function values, we can make $\frac{N-1}{2} \leq k \leq N - 1$ iterations, after which we will have an interval of uncertainty of length $\frac{b-a}{2^k}$.

Since we take the middle of the last interval as the minimum point, the error does not exceed half of the last interval of uncertainty, we get: $\varepsilon_D \leq \frac{b-a}{2\cdot 2^k} \leq \frac{b-a}{2\cdot 2^{\frac{N-1}{2}}}$.

Similarly, in the trichotomy method, we calculate no more than five values of $f(x)$ in the first step and no more than four in the following steps.

So $k$ iterations require $N \leq 5 + 4(k-1) = 4k+1$ ($N \geq 3 + 2(k-1) = 2k+1$) calculations of the $f(x)$. Using $N$ calculations of the function values, we can do $\frac{N-1}{4} \leq k \leq \frac{N-1}{2}$ iterations and achieve the accuracy $\varepsilon_T \leq \frac{b-a}{2\cdot 3^{\frac{N-1}{4}}}$.

**Note 2:**

Since $2^{\frac{N-1}{2}} > 3^{\frac{N-1}{4}} = \sqrt{3}^{\frac{N-1}{2}}$, we have $\frac{b-a}{2\cdot 2^{\frac{N-1}{2}}} \leq \frac{b-a}{2\cdot 3^{\frac{N-1}{4}}}$. As we can see, the estimate for the accuracy given by the interval halving method is lower (better) than the estimate for the trichotomy method. Note that we are talking only about estimates, which give the exact values of the error $\varepsilon$ only in the worst case, when the method has to calculate the maximum possible number of function values at each step. In practice, this rarely happens, as evidenced by our numerous experiments, described in the next section.

## V. Numerical simulations

In our first experiment, we compared the number of calculations of the values of the minimized function to determine the minimum point with the given accuracy. The results are given in Table 1. As we can see, in most examples, the trichotomy method requires fewer calculations of the function values (the best values are in red).

**Table 1. Number of function calculations**

| № | $f(x)$ | $[a, b]$ | $\varepsilon$ | $x^* \approx$ | interval halving method $N_D$ | trichotomy method $N_T$ | golden section search $N_Z$ |
|---|---|---|---|---|---|---|---|
| 1 | $e^x + \frac{1}{x}$ | [0,5; 1] | $10^{-3}$ | 0,7 | 15 | 15 | 15 |
| 2 | $\frac{5}{x} + x^2$ | [0,5; 2] | $10^{-6}$ | 1,35(9) | 37 | 31 | 32 |
| 3 | $-\frac{5}{x^2 - 2x + 5}$ | [0,8; 2] | $10^{-7}$ | 1 | 35 | 31 | 36 |
| 4 | $e^{-2x} + \frac{x^2}{2}$ | [0; 1,5] | $10^{-8}$ | 0,6 | 48 | 40 | 42 |
| 5 | $e^{x-1} + \frac{1}{x}$ | [0; 1,5] | $10^{-6}$ | 1 | 31 | 28 | 32 |

| № | f(x) | [a; b] | ε | x* | | | |
|---|---|---|---|---|---|---|---|
| 6 | $x^2 - x*e^{-x}$ | [0; 1] | $10^{-7}$ | 0,28 | 42 | 34 | 36 |
| 7 | $5*x^2 + \dfrac{1}{x}$ | [0; 2,5] | $10^{-5}$ | 0,46 | 31 | 27 | 28 |
| 8 | $e^{-x} + \dfrac{1}{1-x}$ | [−3; 0] | $10^{-6}$ | 0 | 43 | 40 | 33 |
| 9 | $2 - x + x^2$ | [0; 2] | $10^{-8}$ | 0,5 | 53 | 43 | 42 |
| 10 | $-(x*e^{-0.5x})$ | [0; 3] | $10^{-4}$ | 2 | 22 | 20 | 24 |
| 11 | $-(0,2x + \sin 2x)$ | [0; 3] | $10^{-7}$ | 0,84 | 42 | 37 | 38 |
| 12 | $-(\dfrac{1}{x} - e^{-x})$ | [0; 0,5] | $10^{-5}$ | 0 | 16 | 21 | 25 |
| 13 | $e^x + x^2$ | [−1; 0] | $10^{-6}$ | -0,35 | 34 | 27 | 31 |
| 14 | $x^4 + 2x^2 + 4x$ | [−1; 0] | $10^{-4}$ | -0,67(9) | 22 | 20 | 22 |
| 15 | $x^2 + \sin x$ | [−1; 0] | $10^{-8}$ | -0,44(9) | 48 | 41 | 40 |
| 16 | $e^x + \dfrac{1}{x+2}$ | [−1; 1] | $10^{-5}$ | -0,63 | 30 | 25 | 28 |
| 17 | $\dfrac{-x + (x+2)}{x^2}$ | [−2; 0] | $10^{-8}$ | -2 | 28 | 35 | 42 |
| 18 | $-(5x^2 * e^{-0.5x})$ | [2; 6] | $10^{-7}$ | 4 | 51 | 33 | 39 |
| 19 | $-(0.1x + \cos x)$ | [4; 9] | $10^{-5}$ | 6,38 | 34 | 26 | 30 |
| 20 | $-(\dfrac{\cos 1.5x}{\sin 1.5x} - x^2)$ | [4; 9] | $10^{-6}$ | 4,19 | 31 | 29 | 35 |

In our second experiment, we considered a number of problems for finding a minimum of a unimodal function on an interval in which we know the exact value of the minimum point in advance. In this series of experiments, we limited the methods to the same number N of calculations of the function values and compared the accuracy they achieved to determine the minimum point.

Since the Fibonacci method is the most efficient method in this conditions and provides a maximum reduction of the interval of uncertainty in the given number of trials of function evaluations [2], we compared our method (and the interval halving method) with it. The results are presented in Table 2 and suggest that the trichotomy method is more effective.

Note that Table 2 does not show the last interval of uncertainty, but actually its middles, so there is no contradiction with the fact that the Fibonacci method gives the best reduction of the interval of uncertainty.

**Table 2. The accuracy of the minimum point**

| № | Examples | | | | interval halving method | trichotomy method | Fibonacci search |
|---|---|---|---|---|---|---|---|
| | $f(x)$ | $[a, b]$ | $x^*$ | $N$ | $\varepsilon_D$ | $\varepsilon_T$ | $\varepsilon_F$ |
| 1 | $(x - 1.1)^2$ | $[0; 2]$ | 1.1 | 10 | $0{,}625 * 10^{-2}$ | $0{,}123 * 10^{-2}$ | $0{,}511 * 10^{-2}$ |
| | | | | 20 | $0{,}977 * 10^{-4}$ | $0{,}152 * 10^{-4}$ | $0{,}365 * 10^{-4}$ |
| | | | | 30 | $0{,}611 * 10^{-5}$ | $0{,}019 * 10^{-5}$ | $0{,}082 * 10^{-5}$ |
| 2 | $-(5x^2 \cdot e^{-0.5x})$ | $[1; 6]$ | 4 | 10 | $0{,}313 * 10^{-1}$ | $0{,}062 * 10^{-1}$ | $0{,}056 * 10^{-1}$ |
| | | | | 20 | $0{,}488 * 10^{-3}$ | $0{,}076 * 10^{-3}$ | $0{,}411 * 10^{-3}$ |
| | | | | 30 | $0{,}763 * 10^{-5}$ | $0{,}094 * 10^{-5}$ | $0{,}037 * 10^{-5}$ |
| 3 | $\cos x$ | $[2; 4]$ | $\pi$ | 10 | $0{,}146 * 10^{-1}$ | $0{,}058 * 10^{-1}$ | $0{,}018 * 10^{-1}$ |
| | | | | 20 | $0{,}009 * 10^{-3}$ | $0{,}154 * 10^{-3}$ | $0{,}029 * 10^{-3}$ |
| | | | | 30 | $0{,}635 * 10^{-5}$ | $0{,}027 * 10^{-5}$ | $0{,}015 * 10^{-5}$ |

## VI. Conclusions

This paper proposes a new method for the one-dimensional optimization – the trichotomy method. The new algorithm is a modification of the interval halving method and is based on dividing the interval by five points into six equal parts. The trichotomy method provides the reducing of the interval of uncertainty exactly three times in every iteration. The method's performance has been tested and compared to other algorithms: the interval halving method, golden section method, Fibonacci search. The results show that the trichotomy method is superior to these methods. It usually require less calculations of the values of the minimized function to determine the minimum point with the given accuracy. Also it has better accuracy when finding the minimum point after using *N* calculations of the function values.